\documentclass[aip,jmp,amsmath,amssymb,reprint]{revtex4-1}

\usepackage{amsthm}
\usepackage{graphicx}
\usepackage{tabularx}
\usepackage{latexsym}
\usepackage{enumerate}
\usepackage{times}
\usepackage[T1]{fontenc}

\usepackage{booktabs} 
\usepackage{xfrac}    
\usepackage{hyperref}

\newtheorem{teor}{Theorem}[section]

\theoremstyle{definition}
\newtheorem{defi}{Definition}[section]

\begin{document}
	\title{Scalar conformal invariants of weight zero}

	
	\author{Ignacio S\'anchez-Rodr\'{\i}guez}
	\email{ignacios@ugr.es}
	\homepage[Web:\ ]{http://www.ugr.es/~ignacios}
	\affiliation{Departamento de Geometr\'{\i}a
		y Topolog\'{\i}a, Universidad de Granada,\\ E-18071, Granada, Spain}

	\date{\today}

\begin{abstract}
In the class of metrics of a generic conformal structure
there exists a distinguishing metric.
This was noticed by Albert Einstein in a lesser-known
paper \cite{Eins1921} of 1921. 
We explore this finding from a geometrical
point of view. Then, we obtain
a family of scalar conformal invariants of
weight 0 for generic pseudo-Riemannian conformal structures $[g]$
in more than three dimensions.
In particular, we define
the conformal scalar curvature of $[g]$
and calculate it for some well-known conformal spacetimes,
comparing the results with the Ricci scalar and the Kretschmann scalar.
In the cited paper, Einstein also announced that it is possible to add
an scalar equation to the field equations of General Relativity.
\end{abstract}


\maketitle

\section{Introduction\label{intro}}

We are accustomed to think that in a 
conformal pseudo-Riemannian class of metrics
there is no preferred metric.
\emph{This is not true} for reasonably generic 
conformal structures in four or more dimensions
(neither is it true for three dimensions \cite{Kru1}). 	

In a lesser-known
paper in March 1921 Einstein wrote
\cite{Eins1921} (my notes in squared brackets):

\begin{quote}
	If we put $g^{\prime}_{\mu\nu}=J g_{\mu\nu}$
	[$J$ is given next]
	then $d\sigma^2=g^{\prime}_{\mu\nu}dx_\mu dx_\nu$
	is an invariant that depends only upon
	the ratios of the $g_{\mu\nu}$
	[{the conformal structure }$[g]\,$].
	All Riemann tensors formed as fundamental invariants
	from $d\sigma$ in the customary manner are---when 
	seen as functions of the $g_{\mu\nu}$---Weyl 
	tensors of weight 0. [\dots] Therefore,
	to every law of nature $T(g)=0$ of the
	general theory of relativity,
	there corresponds a law $T(g^{\prime})=0$, 
	which contains only the
	ratios of the $g_{\mu\nu}$.
\end{quote} 

We are going to develop the idea
of Einstein and we will expose a collection
of scalar conformal invariants, which, to my knowledge, 
has been no discussed up to now.
We take the following steps:
\begin{itemize}
	\item Describe the concept of conformal invariant that
	we are going to consider.
	\item Emphasize the existence of a distinguishing metric
	$g^\prime $ in the class of a generic conformal structure $[g]$.
	\item Establish that scalar invariants of 
	$g^\prime $ of $r$-order
	are scalar conformal invariants of $[g]$ of $(r+2)$-order
	\item In particular, I show that the Ricci scalar
	of $g^\prime $ of $[g]$ is a fourth order 
	scalar conformal invariant of $[g]$, which we called 
	\emph{conformal scalar curvature}.
	Then, for some well-known spacetimes, we compare
	the results with some usual scalar metric invariants.
\end{itemize}

\section{Bundle of $G$-structures\label{secBun}}
Given a closed subgroup $G$ of 
$Gl_n\equiv GL(n,\mathbb{R})$,
we get the \emph{bundle $M_G$ of $G$-structures}
on a manifold $M$ that is an associated bundle
to the linear frame bundle $LM\equiv F^1M$, 
which typical fiber is $\sfrac{Gl_n}{G}$.
Local sections of $M_G$ 
correspond to local $G$-structures on $M$.

The \emph{bundle of metrics} on $M$ is 
obtained by taking $G=O_n$, 
the ortho\-gonal group
of a given signature $(p,q)$.
There is a bijective correspondence 
between $\sfrac{Gl_n}{O_n}$ 
and the set of symmetric invertible
matrices of signature $(p,q)$;
therefore, when we give a chart, we can recognize a
pseudo-Riemannian metric $g$ of $M$ as
a section of $M_{O_n}$ given
by a matrix $(g_{ij})$
over the domain of the chart.
We write $\sigma_g\colon M
\to M_{O_n}$ for the section
corresponding to the metric $g$.
 
The \emph{bundle of conformal structures}
arises when $G=C_n:=\mathbb{R}^+\cdot O_n$, 
the conformal ortho\-gonal group
of signature $(p,q)$. 
There is a bijective correspondence 
between $\sfrac{Gl_n}{C_n}$
and the set of symmetric invertible
matrices of signature $(p,q)$
and absolute value of its determinant one.
Now, we can recognize 
a conformal structure $[g]$ on $M$
as a section of $M_{C_n}$,
which is given by the matrix
$(c_{ij})=|\det (g_{ij})|^{-\sfrac{1}{n}}(g_{ij})$ 
over the domain of a chart.
We write $\sigma_{[g]}\colon M
\to M_{C_n}$ for the section
corresponding to the conformal
structure $[g]$.

Let $F^rM$ be the frame bundle of $r$-th order---an 
$r$-frame is an $r$-jet at 0 of 
inverses of charts of $M$---that 
is a principal bundle with structural group
the \emph{jet group} $G_n^r$---group 
of $r$-jets at 0 of diffeomorphisms
of $\mathbb{R}^n$ fixing 0. 
The bundle $M_G^r$ of $r$-jets
of local sections of $M_G$ 
(i.e., $r$-jets of local $G$-structures on $M$)
is an associated bundle to the frame bundle $F^{r+1}M$,
which typical fiber is the space
$(\mathbb{R}^n_G)^r_0\equiv  
J^r_{0}(\mathbb{R}^n,\sfrac{Gl_n}{G})$
of $r$-jets at $0$ of maps
of $\mathbb{R}^n$ to $\sfrac{Gl_n}{G}$
(i.e., $r$-jets of $G$-structures 
on a neighborhood of $0$).

\section{Invariants of $G$-structures\label{secInv}}

A paradigm of scalar invariant for a differential 
geometric structure is the scalar
curvature, or Ricci scalar, $R_g\colon M\to \mathbb{R}$ 
of a Riemannian manifold
$(M, g)$. For $m\in M$, $R_g(m)$ is made from
the partial derivatives up to second order of the metric
$g$ at $m$ and does not depend on which chart has
been used to perform the derivatives, in other words,
$R_g(m)$ is a function of the 2-jet of $g$ at $m$.
Therefore, the Ricci scalar is defined by the function
$R\colon M^2_{O_n}\to \mathbb{R}$,
$R(j^2_m\sigma _g)=R_g(m)$,
with $\sigma _g$ being the section of
$M_{O_n}$ corresponding to the metric $g$.

The Ricci scalar of any metric $g$ 
is said \emph{invariant by
diffeomorphisms} because, for any
local diffeomorphism $\varphi$ of $M$,
it is verified that
$R_{\varphi ^*g}=R_g\circ\varphi $.
This property is equivalent to say 
\cite{ApDifInv} that
$R\circ \widehat{\varphi }^{2}=R$,
$\forall \varphi $,
with $\widehat{\varphi }^{2}$
being the case $r=2$ of
the typical action of a diffeomorphism
$\varphi $ on $M_G^r$, which is defined by:
\begin{equation}
\widehat{\varphi }^{r}\colon
M_G^r \to M_G^r,
\quad
j^r_p\sigma \mapsto
j^r_{\varphi (p)}(\bar{\varphi }
\circ \sigma \circ \varphi ^{-1}),
\end{equation}
with $\bar{\varphi }$ being the standard lifting
of $\varphi $ to $M_G$.	

Following the model of the Ricci scalar function 
$R\colon M^2_{O_n}\to \mathbb{R}$, we define:

\begin{defi}\label{defSca} 
An \emph{scalar invariant 
of $r$-order of $G$-structures on $M$}
is a function
$f\colon M_G^r\to \mathbb{R}$
verifying $f\circ \widehat{\varphi }^r=f$
over  the domain of $\varphi $, for all local
diffeomorphism $\varphi $ of $M$.
\end{defi}

\section{Main results\label{secMai}}
In the setting of conformal differential geometry
there is another concept of invariance,
Let $M$ be a manifold with a conformal structure
$[g]$. A tensor $T$ on $M$,
obtained from the $r$-jet of $g$ by an 
specific formula, is said 
\emph{conformally invariant}---or 
invariant by conformal transformations---if, 
$\forall \bar g\in [g]$, $T$ 
is equal to the tensor $\bar T$
obtained from the $r$-jet of  
$\bar g$ by the same formula.

It is well known that the conformal curvature 
or Weyl tensor $C^i_{\;jkl}$
is conformally invariant.
We define the \emph{square of the Weyl tensor}
to be the scalar $H=C_{ijkl}C^{ijkl}$, 
which is not conformally invariant
because it depends of the metric 
$g$ used to raising and lowering indices.
In this context, a point $p\in M$ 
is called \emph{generic} if $H(p)\neq 0$ for some 
(and then for all) $g\in [g]$; and
a conformal structure $[g]$ is called
\emph{generic}
if all the points of $M$ are generic.

The following result is obtained by Einstein \cite{Eins1921}
with the help, as he himself says, of the mathematician 
Wilhelm Wirtinger.
\begin{teor}\label{teoEin}
Let $[g]$ be a generic 
conformal structure on $M$. 
The metric $g^\prime :=J g\in [g]$,
where $J =\left|H\right|^{1/2}$,
do not depends of the metric $g$
in the class of $[g]$.
\end{teor}

In other words, in the class 
of metrics of a generic 
conformal structure 
\emph{there exists a preferred metric}. 
The theorem 
can be formulated
in a neighborhood of a 
generic point.

We need to calculate the dependence on $g$ 
of the factor $J=|H|^{\sfrac{1}{2}}$.
We can write:
$$H=C_{ijkl}C^{ijkl}=g_{ia}C^a_{\;jkl}
C^i_{\;bcd}g^{bj}g^{ck}g^{dl}.$$
Now, if we change the metric $g$ for $\bar g=\alpha g$ 
then $\bar g_{ij}=\alpha g_{ij}$ and
$\bar g^{ij}=\alpha^{-1} g^{ij}$, and we obtain:
\begin{equation}\label{eqH}
\bar H=\bar C_{ijkl}\bar C^{ijkl}=
\bar g_{ia}C^a_{\;jkl}C^i_{\;bcd
}\bar g^{bj}\bar g^{ck}\bar g^{dl}
=\alpha^{-2}H,
\end{equation}
because of the Weyl tensor is 
conformally invariant:
$\bar C^i_{\;jkl}=C^i_{\;jkl}$.
Then, taking absolute values and 
square roots in Eq.~(\ref{eqH}), 
we obtain $\bar J=\alpha^{-1}J$. 
Finally, we get the result:
\begin{equation}\label{eqJ}
g^\prime=Jg=\alpha^{-1}J\alpha g=\bar J\bar g.
\end{equation}

Furthermore, the scalar $J^\prime$ corresponding to 
$g^\prime$ is constant one, hence we can say that
\emph{the preferred metric of a generic conformal structure
is the only metric in the class that normalizes to $\pm 1$
the square of the Weyl tensor}.

\begin{teor}\label{teoInv} To each scalar metric 
invariant of $r$-order on $M$ corresponds an 
scalar conformal invariant of $(r+2)$-order on $M$.
\end{teor}

In fact, an scalar invariant of $r$-order 
of the preferred metric $g^\prime \in [g]$
is a conformal invariant of $[g]$
of $(r+2)$-order. In particular, the Ricci scalar
of the metric $g^\prime $ is an scalar conformal 
invariant of fourth order that we call the 
\textit{conformal scalar curvature} of $[g]$.

We will return to this theorem 
at the end of the section,
but first let us go back to the notion of 
invariant of Def.~\ref{defSca}.
The action of the group $G_n^{r+1}$
on $(\mathbb{R}^n_G)^r_0$,
which makes $M_G^r$ into an associated bundle 
to $F^{r+1}M$, is defined \cite{ApDifInv} by
\begin{equation}
j^{r+1}_{_0}\xi \cdot j^r_{_0}\mu :=
j^r_{_0}\left( (D\xi \cdot \mu)
\circ \xi ^{-1}\right)
\end{equation}
being $\xi $ a local diffeomorphism of 
$\mathbb{R}^n$ with $\xi (0)=0$,
$\mu $ a smooth map from a neighborhood
of $0$ to $\sfrac{Gl_n}{G}$ (for $G$ closed),
and the last dot for the action of
$Gl_n$ on $\sfrac{Gl_n}{G}$.
Its orbit space is
$\sfrac{(\mathbb{R}^n_G)^r_0}{G^{r+1}_n}$,
consisting of classes of $r$-jets at 0 of 
$G$-structures over $\mathbb{R}^n$
``modulo diffeomorphisms''. 

It is proved \cite{ApDifInv} that 
$r$-order scalar invariants 
of $G$-structures on $M$
are in a natural bijective correspondence
with functions
$\psi\colon \sfrac{(\mathbb{R}^n_G)^r_0}{G^{r+1}_n}
\to \mathbb{R}$
such that $\psi\circ \pi $ is smooth,
with $\pi $ being the projection
of the quotient space. 
We will call a such $\psi$  
an \emph{scalar $G$-invariant of $r$-order}
(it does not depend on $M$!).
We are stating that the nature of an invariant
of $r$-order can be appreciated when 
it works on $r$-jets at 0
of $G$-structures on $\mathbb{R}^n$
in a neighborhood of 0. 

With this simplified concept of scalar invariant
it is easier to prove the following
properties of invariants (in the rest of the section 
all the invariants are scalar invariants):
\begin{itemize}
	\item\emph{ An invariant of $r$-order is 
	invariant of $(r+s)$-order}, $\forall s$.
	Because there exists the projection
	\[\Pi\colon (\mathbb{R}^n_G)^{r+s}_0
	\to (\mathbb{R}^n_G)^r_0,\quad 
	j^{r+s}_0\mu\mapsto j^r_0\mu,\]
	which passes to a map $\Pi^\prime$ 
	between quotient spaces, we obtain
	the result on invariants from:
	\[\sfrac{(\mathbb{R}^n_G)^{r+s}_0}{G^{r+s+1}_n}
	\stackrel{\Pi^\prime}{\longrightarrow}
	\sfrac{(\mathbb{R}^n_G)^r_0}{G^{r+1}_n}
	\stackrel{\psi}{\longrightarrow}\mathbb{R}.\]
	
	\item \emph{If $U\subset G$ is a subgroup,
	a $G$-invariant of $r$-order is a 
	$U$-invariant of $r$-order}. In effect: the map 
	$\beta \colon \sfrac{Gl_n}{U}
	\to \sfrac{Gl_n}{G} ,\  aU\mapsto aG$, 
	verifies $\beta \circ L_b=L_b\circ \beta$, 
	$\forall b\in Gl_n$; 
	hence the map
	\[\gamma\colon (\mathbb{R}^n_U)^r_0
	\to (\mathbb{R}^n_G)^r_0,\quad 
	j^r_0\mu\mapsto j^r_0(\beta\circ\mu)\]
	passes to a map $\gamma^\prime$ 
	between quotient spaces, and we obtain
	this result on invariants from:
	\[  \sfrac{(\mathbb{R}^n_U)^r_0}{G^{r+1}_n}
	\stackrel{\gamma^\prime}{\longrightarrow}
	\sfrac{(\mathbb{R}^n_G)^r_0}{G^{r+1}_n}
	\stackrel{\psi}{\longrightarrow}\mathbb{R}.\]
	Therefore, \emph{conformal invariants of a 
	given order are metric invariants of the same order}.
	Moreover, the conformal invariants
	are the metric invariants which are
	conformally invariant.
	
	\item \emph{A metric invariant of $r$-order is a 
	conformal invariant of $(r+2)$-order}
	(hence Theorem \ref{teoInv} holds).
	Because $J$ is made from
	the partial derivatives 
	up to second order of the metric
	$g$ of $M$, the well-defined map 
	between conformal structures and
	metrics, given by
	$[g]\mapsto g^\prime=Jg$, depends
	of the 2-jet of the components
	$g_{ij}$ at each point of $M$ 
	or, equivalently, of the 2-jet of 
	the components $c_{ij}$ that characterize $[g]$. 
	The same is true for conformal
	structures $[g]$ on $\mathbb{R}^n$
	in a neighborhood of 0; then
	the map between conformal structures
	and metrics
	induces, for successive orders, the maps
	\[\Gamma\colon (\mathbb{R}^n_{C_n})^{r+2}_0
	\to (\mathbb{R}^n_{O_n})^r_0,
	\quad j^{r+2}_0\sigma_{[g]}\mapsto j^r_0 \sigma_{g^\prime},\]
	each one of which passes to a map $\Gamma^\prime$ 
	between quotient spaces, and we obtain
	the result from:
	\[\sfrac{(\mathbb{R}^n_{C_n})^{r+2}_0}{G^{r+3}_n}
	\stackrel{\Gamma^\prime}{\longrightarrow}
	\sfrac{(\mathbb{R}^n_{O_n})^r_0}{G^{r+1}_n}
	\stackrel{\psi}{\longrightarrow}\mathbb{R}.\]
	
\end{itemize}

\section{Conformal scalar curvature\label{secCSC}}

I have calculated the Ricci scalar
of $g^\prime$, which we have called
\emph{conformal scalar curvature
of} $[g]$ and denoted by $S_{[g]}$, 
for some well known 
spacetime metrics $g$. To do this 
for a metric
$g$, we calculate $J$ and
the Ricci scalar of
$g^\prime=J g$. I used
the package  \emph{xAct\`{}xCoba\`{}} 
for Mathematica.
In the Table \ref{tab1} below, 
I give some results
comparing the Ricci scalar 
and Kretschmann scalar of $g$
with the conformal scalar curvature 
of $[g]$;
this table shows the calculations for
the following metrics:
\begin{itemize}
	\item Schwarzschild: $-(1-\frac{2M}{r})
	dt^2+(1-\frac{2M}{r})^{-1}dr^2
	+r^2(d\theta ^2+\sin ^2\theta d\phi ^2)$
	\item Reissner-Nordstr\o m: 
	$-(1-\frac{2M}{r}+\frac{q^2}{r^2})dt^2
	+(1-\frac{2M}{r}+\frac{q^2}{r^2})^{-1}dr^2
	+r^2(d\theta ^2+\sin ^2\theta d\phi ^2)$
	\item G\"{o}del: $-dt^2+
	(1+\frac{r^2}{4a^2})^{-1}dr^2
	+r^2(1-\frac{r^2}{4a^2})d\phi ^2
	-\frac{\sqrt{2}\,r^2}{a}\,dt\,d\phi +dz^2$
	\item Barriola-Vilenkin: $-dt^2+dr^2
	+k^2r^2(d\theta ^2+\sin ^2\theta d\phi ^2)$
\end{itemize}

\begin{table}[h]
	\caption{Ricci scalar $R_g$,
	Kretschmann scalar $K_g$
	and conformal scalar curvature 
	$S_{[g]}$ for some solutions $g$ of the
	Einstein equations.\label{tab1}}

\begin{ruledtabular}
		\begin{tabularx}{\textwidth}{@{\extracolsep{\fill}}lccc}
		\emph{Metrics} & $R_g$ & $K_g$ & $S_{[g]}$\\
		\midrule
		Schwarzschild & 0 & $\frac{48M^2}{r^6}$ & 
		$\frac{9\sqrt{3}}{4}(1-\frac{r}{6M})$ \\[2mm]
		Reissner-Nordstr\o m & 0 & $\frac{8(7q^4-12Mq^2r
		+6M^2 r^2)}{r^8}$ & $\frac{9\sqrt{3}}{8}(1-
		\frac{r}{3M}+\frac{q^2}{3M^2})$
		\\[2mm]
		G\"{o}del & $-\frac{1}{a^2}$ & $\frac{3}{a^4}$ 
		& $-\frac{\sqrt{3}}{2}$ 	\\[2mm]
		Barriola-Vilenkin  & $\frac{2(1-k^2)}{k^2r^2}$ &
		$\frac{4(1-k^2)^2}{k^4r^4}$ & $-\sqrt{3}$\\[2mm]
	\end{tabularx}
\end{ruledtabular}
\end{table}

Let me highlight that replacing the metric $g$ 
by a conformally related metric $\alpha g$, 
Ricci and Kretschmann scalars
get complicated with terms including
up to second order derivatives of $\alpha $.
In contrast, the conformal scalar curvature
keeps unalterable, meaning that its weight is zero
in the Weyl's sense \cite{FefGra}.

\section{Conclusion}
At the end of his article \cite{Eins1921}, 
Einstein proposed the obvious 
addition of the differential equation 
$J =J_{_0}$ (constant), 
free of genericity conditions, 
to the field equations of General Relativity
because it is 
\textquotedblleft\emph{a logical possibility 
that is worthy of publication, 
which may be useful for physics 
or not.}\textquotedblright

In dimension four, the existence of a distinguishing metric
in the class of a conformal structure with $H\neq 0$ arises also
from the theory of Weyl gravity (or conformal gravity).
This theory is governed by the Lagrangian
$H\,\Omega_g$, with $\Omega_g$ being the metric volume form.
This Lagrangian, in dimension 4, is invariant by
conformal transformations
(in dimension greater than 4 
what is conformally invariant is 
$|H|^{\sfrac{n}{4}}\,\Omega_g$,
which is the former Lagrangian, except sign, when $n=4$).
Then $H \Omega_g$ 
is a distinctive volume form provided only
by the conformal structure $[g]$.
Therefore, it can be concluded that
there exists a unique metric
$g^\prime \in [g]$ such that $\Omega_{g^\prime}=H \Omega_g$.

These previous facts suggest new perspectives of research
using that preferred metric of the conformal class,
which is surprisingly almost unknown in the research 
bibliography (I only found one 
reference---not on subject of history---to this paper of Einstein 
pointing at this distinguishing metric \cite{Tred1}).

\begin{acknowledgments}
This work was presented for the first time at the 
Spanish-Portuguese Relativity Meeting -- EREP 2017 
held in M\'alaga (Spain), 12-15 September 2017. 
I am very grateful to all the colleagues at the meeting who,
through their comments, have contributed to improve this paper.
This research was partially supported by 
Junta de Andaluc\'{\i}a (Spain),
under grant P.A.I. FQM-324.
\end{acknowledgments}

\bibliography{ConfInvMetric}

\providecommand{\noopsort}[1]{}\providecommand{\singleletter}[1]{#1}%
\begin{thebibliography}{5}%
\makeatletter
\providecommand \@ifxundefined [1]{%
 \@ifx{#1\undefined}
}%
\providecommand \@ifnum [1]{%
 \ifnum #1\expandafter \@firstoftwo
 \else \expandafter \@secondoftwo
 \fi
}%
\providecommand \@ifx [1]{%
 \ifx #1\expandafter \@firstoftwo
 \else \expandafter \@secondoftwo
 \fi
}%
\providecommand \natexlab [1]{#1}%
\providecommand \enquote  [1]{``#1''}%
\providecommand \bibnamefont  [1]{#1}%
\providecommand \bibfnamefont [1]{#1}%
\providecommand \citenamefont [1]{#1}%
\providecommand \href@noop [0]{\@secondoftwo}%
\providecommand \href [0]{\begingroup \@sanitize@url \@href}%
\providecommand \@href[1]{\@@startlink{#1}\@@href}%
\providecommand \@@href[1]{\endgroup#1\@@endlink}%
\providecommand \@sanitize@url [0]{\catcode `\\12\catcode `\$12\catcode
  `\&12\catcode `\#12\catcode `\^12\catcode `\_12\catcode `\%12\relax}%
\providecommand \@@startlink[1]{}%
\providecommand \@@endlink[0]{}%
\providecommand \url  [0]{\begingroup\@sanitize@url \@url }%
\providecommand \@url [1]{\endgroup\@href {#1}{\urlprefix }}%
\providecommand \urlprefix  [0]{URL }%
\providecommand \Eprint [0]{\href }%
\providecommand \doibase [0]{http://dx.doi.org/}%
\providecommand \selectlanguage [0]{\@gobble}%
\providecommand \bibinfo  [0]{\@secondoftwo}%
\providecommand \bibfield  [0]{\@secondoftwo}%
\providecommand \translation [1]{[#1]}%
\providecommand \BibitemOpen [0]{}%
\providecommand \bibitemStop [0]{}%
\providecommand \bibitemNoStop [0]{.\EOS\space}%
\providecommand \EOS [0]{\spacefactor3000\relax}%
\providecommand \BibitemShut  [1]{\csname bibitem#1\endcsname}%
\let\auto@bib@innerbib\@empty
\bibitem [{\citenamefont {Einstein}(1921)}]{Eins1921}%
  \BibitemOpen
  \bibfield  {author} {\bibinfo {author} {\bibfnamefont {A.}~\bibnamefont
  {Einstein}},\ }\bibfield  {title} {\enquote {\bibinfo {title} {{\"{U}}ber
  eine naheliegende erg{\"{a}}nzung des fundamentes der allgemeinen
  relativit{\"{a}}tstheorie},}\ }\href@noop {} {\bibfield  {journal} {\bibinfo
  {journal} {Berl.\ Ber.}\ ,\ \bibinfo {pages} {261--264}} (\bibinfo {year}
  {1921})},\ \bibinfo {note} {``On a Natural Addition to the Foundation of the
  General Theory of Relativity'', \emph{CPAE, Vol. 7: The Berlin Years:
  Writings, 1918-1921 (English translation supplement)}, Princeton University
  Press, Doc. 54}\BibitemShut {NoStop}%
\bibitem [{\citenamefont {Kruglikov}(2017)}]{Kru1}%
  \BibitemOpen
  \bibfield  {author} {\bibinfo {author} {\bibfnamefont {B.}~\bibnamefont
  {Kruglikov}},\ }\bibfield  {title} {\enquote {\bibinfo {title} {Conformal
  differential invariants},}\ }\href@noop {} {\bibfield  {journal} {\bibinfo
  {journal} {J.\ Geom.\ Phys.}\ }\textbf {\bibinfo {volume} {113}},\ \bibinfo
  {pages} {170--175} (\bibinfo {year} {2017})}\BibitemShut {NoStop}%
\bibitem [{\citenamefont {S{\'{a}}nchez-Rodr{\'{\i}}guez}(2017)}]{ApDifInv}%
  \BibitemOpen
  \bibfield  {author} {\bibinfo {author} {\bibfnamefont {I.}~\bibnamefont
  {S{\'{a}}nchez-Rodr{\'{\i}}guez}},\ }\href@noop {} {\enquote {\bibinfo
  {title} {An approach to differential invariants of {$G$}-structures},}\ }
  (\bibinfo {year} {2017}),\ \Eprint {http://arxiv.org/abs/1709.02382}
  {arXiv:1709.02382 [math.DG]} \BibitemShut {NoStop}%
\bibitem [{\citenamefont {Fefferman}\ and\ \citenamefont
  {Graham}(2012)}]{FefGra}%
  \BibitemOpen
  \bibfield  {author} {\bibinfo {author} {\bibfnamefont {C.}~\bibnamefont
  {Fefferman}}\ and\ \bibinfo {author} {\bibfnamefont {C.~R.}\ \bibnamefont
  {Graham}},\ }\href@noop {} {\emph {\bibinfo {title} {The Ambient Metric
  (AM-178)}}}\ (\bibinfo  {publisher} {Princeton University Press},\ \bibinfo
  {year} {2012})\BibitemShut {NoStop}%
\bibitem [{\citenamefont {Treder}(1992)}]{Tred1}%
  \BibitemOpen
  \bibfield  {author} {\bibinfo {author} {\bibfnamefont {H.~J.}\ \bibnamefont
  {Treder}},\ }\bibfield  {title} {\enquote {\bibinfo {title} {Conform
  invariance and {M}ach\'{}s principle in cosmology},}\ }\href@noop {}
  {\bibfield  {journal} {\bibinfo  {journal} {Found.\ Phys.}\ }\textbf
  {\bibinfo {volume} {22}},\ \bibinfo {pages} {1089--1093} (\bibinfo {year}
  {1992})}\BibitemShut {NoStop}%
\end{thebibliography}%

\end{document}